\documentclass[letterpaper,conference]{ieeeconf}

\IEEEoverridecommandlockouts
\usepackage{arydshln}
\usepackage{graphicx}
\usepackage{amsmath}
\usepackage{amssymb, verbatim}
\usepackage{mathrsfs}
\usepackage{dsfont}
\usepackage{url}
\usepackage[usenames,dvipsnames]{color}
\usepackage[compress]{cite}
\usepackage{mdwlist}
\usepackage{paralist}
\usepackage{algorithmic}
\usepackage{algorithm}
\usepackage{cleveref}
\usepackage{stmaryrd}
\usepackage{tikz}
\usepackage{tkz-graph}
\usepackage{lipsum}

\usetikzlibrary{arrows,positioning,automata,calc}
\usetikzlibrary{shapes,snakes,patterns}

\newtheorem{lemma}{\sc Lemma}
\newtheorem{theorem}[lemma]{\sc Theorem}

\newtheorem{remark}{\sc Remark}
\newtheorem{assumption}{\sc Assumption}

\renewcommand{\matrix}[2]{\left[\begin{array}{#1} #2 \end{array}\right]}

\DeclareMathOperator*{\argmin}{arg\,min}

\DeclareMathOperator*{\diag}{diag}
\DeclareMathOperator*{\trace}{trace}

\newcommand{\q}{\llbracket q\rrbracket}
\newcommand{\n}[1]{\llbracket n_{#1}\rrbracket}

\floatname{algorithm}{Algorithm}

\renewcommand{\footnoterule}{%
  \kern -7pt
  \hrule width 0.3\textwidth height .5pt
  \kern 2pt
}

\renewcommand{\baselinestretch}{0.96} 

\begin{document}

\title{SiMpLIfy: A Toolbox for Structured Model Reduction}
\author{
Martin Biel, Farhad~Farokhi, and Henrik~Sandberg
\thanks{M.~Biel and H.~Sandberg are with the ACCESS Linnaeus Center, School of Electrical Engineering, KTH Royal Institute of Technology, SE-100 44 Stockholm, Sweden. E-mails: \{mbiel,hsan\}@kth.se}
\thanks{F.~Farokhi is with the Department of Electrical and Electronic Engineering, University of Melbourne, Parkville, Victoria 3010, Australia. E-mail: farhad.farokhi@unimelb.edu.au}
\thanks{The work was supported, in part, by the Swedish Research Council and the Knut and Alice Wallenberg Foundation.} 
}

\maketitle

\begin{abstract} In this paper, we present a toolbox for structured model reduction developed for MATLAB. In addition to structured model reduction methods using balanced realizations of the subsystems, we  introduce a numerical algorithm for structured model reduction using a subgradient optimization algorithm. We  briefly present the syntax for the toolbox and its features. Finally, we demonstrate the applicability of various model reduction methods in the toolbox on a structured mass-spring mechanical system.
\end{abstract}

\section{Introduction}
Recent developments in control engineering and communication networks have enabled us to construct large-scale engineering or physical systems, such as smart grids and intelligent transportation systems, that are intertwined with our daily life. These systems are, most often, composed of several smaller units that are interconnected to each other through their dynamics, controllers, or performance criteria. The interconnection patterns are, typically, governed by the physical characteristics of the system and the geographical distribution of its subsystems. Since, most often, these systems are scattered across vast areas, the interconnection pattern is structured and carries valuable insights about the weaknesses and the strengths of the system. Because of the large scale of these systems, it is desirable to develop model (order) reduction methods that can reduce the order of the system while preserving the original interconnection between the subsystems (to not sacrifice the mentioned insights). 

Model reduction has been extensively studied in the past~\cite{moore1981principal,glover1984all,obinata2001model}, however, most of these methods do not preserve the interconnection structure between the subsystems. Therefore, several studies have specifically focused on structured model reduction. An early study in~\cite{enns1984model} proposes a heuristic method for frequency-weighted model reduction. Note that the frequency-weighted model reduction can be seen as a structured model reduction for serial interconnection patterns. The idea was later generalized to feedback interconnection pattern in~\cite{schelfhout1996note}. These ideas were generalized to  structured model reduction with arbitrary static networks in~\cite{vandendorpe2004model,VandendorpeDooren2008} and with dynamic networks in~\cite{SandbergOCA}. Generalized structured Gramians were introduced using linear matrix inequalities and used for structured model reduction in~\cite{Li2005145,Zhou1995235}. When using generalized structured Gramians, bounds on the reduction error were provided~\cite{Zhou1995235,SandbergOCA}, however, the existence of generalized structured Gramians cannot be guaranteed unless in special cases~\cite{trnka2013structured}. 

In this paper, we present a structured model reduction toolbox for MATLAB. To describe the toolbox properly, we first survey various structured model reduction algorithms in the literature. We start with model reduction methods using the balanced realizations. For extracting the balanced realizations of the subsystems, we use the structured Gramians. The structured Gramians are calculated heuristically by extracting the block-diagonal entries of regular controllability and observability Gramians. We later use the generalized structured Gramians to construct the balanced realizations. The generalized structured Gramians are extracted using linear matrix inequalities. Bounds on the reduction error for balanced truncation using generalized structured Gramians are presented. Considering that these methods do not provide (sub)optimal reduced systems, we present a numerical method using subgradient optimization algorithm. The numerical algorithm builds upon  $\mathcal{H}_\infty$ synthesis results in~\cite{1576856}. After presenting the methods, we present a brief description of the toolbox and its syntax. We also demonstrate the applicability of the included structured model reduction methods on a structured mechanical system.

The rest of the paper is organized as follows. In Section~\ref{sec:formulation}, we present a mathematical framework for presenting interconnected systems and formulate the structured model reduction problem. In Section~\ref{sec:balanced}, we survey the model reduction methods using the balanced realizations of the subsystems. Model reduction using the subgradient optimization algorithm is presented in Section~\ref{sec:subgradient}. Finally, we present the present the numerical example in~Section~\ref{sec:numerical} and conclude the paper in Section~\ref{sec:conclusions}.

\subsection{Notation}
Let $\mathbb{N}$ and $\mathbb{R}$ denote the sets of integers and reals. Furthermore, define $\q=\{1,\dots,q\}$ for any $q\in\mathbb{N}$. Let the Hardy space $\mathcal{H}_\infty$ be the space of complex Lebesgue measurable functions that are analytic and bounded in the right half plane $\{s\in\mathbb{C}\,|\,\mathrm{Re}(s)>0\}$. For all $G(s)\in\mathcal{H}_\infty$, we define the $\mathcal{H}_\infty$-norm as $\|G(s)\|_\infty=\sup_{\omega\in\mathbb{R}} \sigma_{\max} (G(j\omega))$, where $\sigma_{\max}(\cdot)$ denotes the largest singular value of a complex matrix. Moreover, let $\mathcal{RH}_\infty$ be the set of proper rational functions with real coefficients in $\mathcal{H}_\infty$. For a rational transfer function $G(s)$, $\deg(G(s))$ denotes its McMillan degree. Finally, the set of symmetric positive semidefinite matrices in $\mathbb{R}^{n\times n}$ is denoted by $\mathcal{S}_+^n$.

\section{Problem Formulation and Preliminaries} \label{sec:formulation}
\subsection{Interconnected Systems}
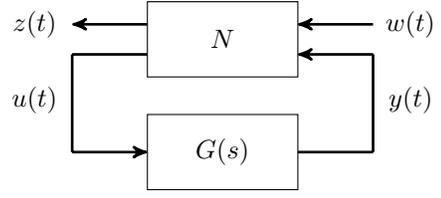
\begin{figure}
\centering
\begin{tikzpicture}[scale=1,>=stealth']
\node (N) at (+0.0,+0.0) [draw,rectangle,fill=white,minimum width=2.0cm,minimum height=1.0cm] {$N$   };
\node (G) at (+0.0,-1.5) [draw,rectangle,fill=white,minimum width=2.0cm,minimum height=1.0cm] {$G(s)$};
\draw[->,line width=1pt] (+2.0,-0.2) to node[right]{} (+2.0,-1.5)
						 (+2.0,-1.5) to node[right]{} (+1.0,-1.5)
						 (+2.0,-0.2) to node[right]{} (+1.0,-0.2);
\draw[->,line width=1pt] (-1.0,-0.2) to node[right]{} (-2.0,-0.2)
						 (-2.0,-0.2) to node[right]{} (-2.0,-1.5)
						 (-2.0,-1.5) to node[right]{} (-1.0,-1.5);
\draw[->,line width=1pt] (+2.0,+0.2) to node[right]{} (+1.0,+0.2);
\draw[->,line width=1pt] (-1.0,+0.2) to node[right]{} (-2.0,+0.2);
\node (1) at (+2.5,+0.2) {$w(t)$};
\node (1) at (-2.5,+0.2) {$z(t)$};
\node (1) at (-2.5,-0.8) {$u(t)$};
\node (1) at (+2.5,-0.8) {$y(t)$};
\end{tikzpicture}
\caption{\label{fig:interconnection} The interconnected systems. The subsystem transfer functions to be reduced $(G_i(s))_{i=1}^n$ are stored in the block-diagonal transfer function $G(s)$. The network structure is captured by the matrix $N$. }
\end{figure}

Here, we present a framework for representing interconnected systems. We follow the convention in~\cite{SandbergOCA} to represent interconnected systems as the feedback form in Figure~\ref{fig:interconnection}. Let $G(s)$ contain the transfer functions of $q\in\mathbb{N}$ subsystems as its block-diagonal entries 
\begin{align*}
G(s)=\diag(G_1(s),\dots, G_q(s)),
%\triangleq \matrix{c|c}{A_G & B_G \\ \hline C_G & D_G},
\end{align*}
where $G_i(s)\in(\mathcal{RH}_\infty)^{p_i\times m_i}$ for $m_i,p_i\in\mathbb{N}$ is the transfer function of subsystem $i\in\q$. For each $i\in\q$, we assume 
$$
G_i(s)=C_i(sI-A_i)^{-1}B_i+D_i,
$$
where $A_i\in\mathbb{R}^{n_i\times n_i}$, $B_i\in\mathbb{R}^{n_i\times m_i}$, $C_i\in\mathbb{R}^{p_i\times n_i}$, and $D_i\in\mathbb{R}^{p_i\times p_i}$ for some $n_i\in\mathbb{N}$. Clearly, we have
$$
G(s)=C_G(sI-A_G)^{-1}B_G+ D_G,
$$
where
\begin{align*}
A_G&=\diag(A_1,\dots,A_q)\in\mathbb{R}^{n\times n},\\
B_G&=\diag(B_1,\dots,B_q)\in\mathbb{R}^{n\times m},\\
C_G&=\diag(C_1,\dots,C_q)\in\mathbb{R}^{p\times n},\\
D_G&=\diag(D_1,\dots,D_q)\in\mathbb{R}^{p\times p},
\end{align*}
with $n=\sum_{i=1}^q n_i$, $m=\sum_{i=1}^q m_i$, and $p=\sum_{i=1}^q p_i$. The augmented input and output vector of subsystems are, respectively, $u(t)\in\mathbb{R}^m$ and $y(t)\in\mathbb{R}^p$. Moreover, let $w(t)\in\mathbb{R}^{m'}$ and $z(t)\in\mathbb{R}^{p'}$ be the external inputs and outputs. The network is modelled using $N\in\mathbb{R}^{(m'+m)\times (p'+p)}$ as
\begin{align*}
\matrix{c}{z(t) \\ y(t)}=N\matrix{c}{w(t) \\ u(t)}
\triangleq \matrix{cc}{ D_E & D_F \\ D_H & D_K }\matrix{c}{w(t) \\ u(t)}.
\end{align*}
Note that the assumption that the network is static is without loss of generality as one can always absorb the network dynamics into the subsystems or introduce additional subsystems. The transfer function matrix of the complete interconnected system is given by the lower linear fractional transformation
\begin{align*}
\mathcal{F}(N,G(s))
&\triangleq D_E+D_F(I-G(s)D_K)^{-1}G(s)D_H(s)\\
&=C(sI-A)^{-1}B+D,
\end{align*}
where
\begin{align*}
A&=A_G+B_G(I-D_KD_G)^{-1}D_KC_G,\\
B&=B_G(I-D_KD_G)^{-1}D_H,\\
C&=D_F(I-D_GD_K)^{-1}C_G,\\
D&=D_E+D_FD_G(I-D_KD_G)^{-1}D_H.
\end{align*}
We make the following assumption throughout the paper.
\begin{assumption} $\mathcal{F}(N,G(s))\in(\mathcal{RH}_\infty)^{p'\times m'}$. \end{assumption}

\begin{remark} In general, model reduction problems are traditionally defined for stable systems due to various reasons. Firstly, the controllability and observability Gramians, that are typically used for balanced realization methods, are not well-defined for unstable systems. Secondly, and more importantly, model reduction of the anti-stable part of the systems (i.e., the part of the system that contains only the unstable poles) is not meaningful as the outputs of two anti-stable systems that do not have the same transfer function drift arbitrarily apart even when excited with the same input. Therefore, we do not address structured model reduction for unstable closed-loop systems. 
\end{remark}

\subsection{Model Reduction Problem}
To extract the reduced system, we propose solving the optimization problem
\begin{align*}
\hat{G}(s)\in\hspace{-.3in}
\argmin_{\footnotesize
\begin{array}{c}
\hat{G}(s)=\diag((\hat{G}_i(s))_{i\in\q}),\\
\hat{G}_i(s)\in(\mathcal{RH}_\infty)^{p_i\times m_i} , \forall i\in\q, \\
\deg(\hat{G}_i(s))\leq r_i, \forall i\in\q
\end{array} }\hspace{-.4in}
\|\mathcal{F}(N,\hat{G}(s))-\mathcal{F}(N,G(s))\|_\infty,
\end{align*}
where $r_i\leq n_i$ is the order of the reduced subsystem~$i\in\q$. Note that this optimization problem is, generally, nonconvex~\cite{dullerud2010course}. Therefore, several heuristic methods have been proposed to find ``good'' solutions. In the next section, we review some of these methods. Subsequently, we propose a numerical algorithm using a subgradient optimization algorithm to find a locally optimal solution. These methods are all implemented in the toolbox.

%[This is a difficult problem; NP hard; cite; we use suboptimal methods that rely on structured Gramians; Then, we introduce a novel method using the subgradients]

\section{Structured Balanced Truncation and Singular Perturbation} \label{sec:balanced}
First, we define the structured Gramians and the generalized structured Gramians. Then, we introduce the balanced truncation and singular perturbation using these Gramians.

\subsection{Structured Gramians} \label{sub:gramians}
In this subsection, we make the following assumption.

\begin{assumption} The state-space representation of the closed-loop system $(A,B,C)$ is minimal.
\end{assumption}

Let the regular controllability Gramian $P'\in\mathbb{R}^{n\times n}$ and the regular observability Gramian $Q'\in\mathbb{R}^{n\times n}$ of the overall system $\mathcal{F}(N,G(s))$ be calculated as the unique positive definite solutions of the Lyapunov equations
\begin{align*}
AP'+P'A^\top+BB^\top=0,\\
A^\top Q'+Q'A+C^\top C=0.
\end{align*}
We can use these regular Gramians to balance the overall system and reduce its dimension~\cite{obinata2001model}, however, doing so, we will lose the inherent structure of the system. Therefore, we need to define structured Gramians. Considering the order of the subsystems, we may decompose the regular Gramians as
\begin{align*}
P'=\matrix{ccc}{
P'_{11} & \cdots & P'_{1q} \\
\vdots & \ddots & \vdots \\
P'_{q1} & \cdots & P'_{qq} }\hspace{-.05in}, \;\;
Q'=\matrix{ccc}{
Q'_{11} & \cdots & Q'_{1q} \\
\vdots & \ddots & \vdots \\
Q'_{q1} & \cdots & Q'_{qq} }\hspace{-.05in},
\end{align*}
where $P'_{ij},Q'_{ij}\in\mathbb{R}^{n_i\times n_j}$ for all $i,j\in\q$. The matrices $P=\diag(P'_{11},\dots,P'_{qq})$ and $Q=\diag(Q'_{11},\dots,Q'_{qq})$ are, respectively, the structured controllability and observability Gramians for the interconnected system.

%The idea of using block-diagonal entries of the Gramian matrices for balancing a system dates back to he results on the frequency-weighted model reduction in~\cite{enns1984model}. This was later generalized to feedback interconnections in~\cite{schelfhout1996note} and, subsequently, to structured model reduction with static networks in~\cite{vandendorpe2004model} and with dynamic networks in~\cite{SandbergOCA}.

\subsection{Generalized Structured Gramians} \label{sub:generalizedgramians}
The generalized structured Gramians (introduced originally in~\cite{beck1996model} for uncertain model reduction) can be extracted from semi-definite programming problems
\begin{align*}
P\in &\argmin_{
\footnotesize
\begin{array}{c}
P=\diag((P_{ii})_{i\in\q})
\end{array} }  && \hspace{-.1in} \trace(P), \\
&\hspace{.53in} \mathrm{s.t.} && \hspace{-.1in} AP+PA^\top+BB^\top\leq 0,\\
&  && \hspace{-.1in} P_{ii}\in\mathcal{S}_{+}^{n_i},\forall i\in\q,
\end{align*}
and
\begin{align*}
Q\in &\argmin_{
\footnotesize
\begin{array}{c}
Q=\diag((Q_{ii})_{i\in\q})
\end{array} }  && \hspace{-.1in} \trace(Q), \\
&\hspace{.53in} \mathrm{s.t.} && \hspace{-.1in} A^\top Q+QA+C^\top C\leq 0,\\
&  && \hspace{-.1in} Q_{ii}\in\mathcal{S}_{+}^{n_i},\forall i\in\q,
\end{align*}
Unfortunately, the generalized structured Gramians may not exist in general unless focusing on specific categories of systems, e.g., the subsystems are strictly positive real~\cite{trnka2013structured}.

\subsection{Balanced Realization}
Now, for each $i\in\q$, we may find transformation $T_i$ so that $T_iP_{ii}T_i^\top=T_i^{-\top}Q_{ii}T_i^{-1}=\Sigma_i\in\mathcal{S}_{+}^{n_i}$, where $P_{ii}$ and $Q_{ii}$ are the block-diagonal entries of either the structured Gramians or the generalized structured Gramians. We assume that the transformation $T_i$ is chosen to guarantee that the diagonal entries of $\Sigma_i$ appear in a descending order. This is without loss of generality as we can always switch the order of the columns in the transformation $T_i$ to achieve such a property. The structured Hankel singular values for subsystem $i\in\q$ can be computed as $\sigma_{i,k}=\sqrt{\lambda_k(P_{ii}Q_{ii})}$ for all $k\in\n{i}$. Clearly, $\Sigma_i=\diag((\sigma_{i,k})_{k\in\n{i}})$. Now, we may find the corresponding state-space representation for this transformation as
\begin{align*}
\bar{A}_i&=T_i^{-1}A_iT_i,&
\bar{B}_i&=T_i^{-1}B_i,&
\bar{C}_i&=C_iT_i,&
\bar{D}_i&=D_i.
\end{align*}
In what follows, we use this balanced realization to extract the reduced subsystems.

\subsection{Balanced Truncation}
We may decompose the model matrices of the balanced subsystem $i\in\q$ as
\begin{align*}
\bar{A}_i=\matrix{cc}{\bar{A}_i^{11} & \bar{A}_i^{12} \\ \bar{A}_i^{21} & \bar{A}_i^{22}}\hspace{-.04in},\,\,
\bar{B}_i=\matrix{cc}{\bar{B}_i^{1} \\ \bar{B}_i^{2}}\hspace{-.04in},\,\,
\bar{C}_i=\matrix{cc}{\bar{C}_i^{1} & \bar{C}_i^{2}}\hspace{-.04in},
\end{align*}
where $\bar{A}_i^{11}\in\mathbb{R}^{r_i\times r_i}$, $\bar{B}_i^{1}\in\mathbb{R}^{r_i\times m_i}$, and $\bar{C}_i^{1}\in\mathbb{R}^{p_i\times r_i}$ with $r_i\in\mathbb{N}$ so that $r_i\leq n_i$. Following this, we can easily calculate the truncated subsystem $i\in\q$ as
$$
\hat{G}_i(s)=\hat{C}_i(sI-\hat{A}_i)^{-1}\hat{B}_i+\hat{D}_i,
$$
where 
\begin{align*}
\hat{A}_i&=\bar{A}_i^{11},&
\hat{B}_i&=\bar{B}_i^{1},&
\hat{C}_i&=\bar{C}_i^{1},&
\hat{D}_i&=\bar{D}_i.
\end{align*}
A nice property of the balanced truncation is that the reduced system and the original system behave similarly for high frequencies, that is, $\mathcal{F}(N,\hat{G}(\infty))=\mathcal{F}(N,G(\infty))$~\cite{green1995linear}.

If we use the structured Gramians in Subsection~\ref{sub:gramians}, we cannot guarantee a good performance or even informative upper bounds on the reduction error (see Theorem~2 in~\cite{SandbergOCA}). However, upon using the generalized structured Gramians in Subsection~\ref{sub:generalizedgramians}, we get the following intuitive bounds on the quality of the reduced model.

\begin{theorem}[\hspace{-.001in}\cite{SandbergOCA}] Let $\hat{G}(s)=\diag((\hat{G}_i)_{i\in\q})$, where $\hat{G}_i(s)$, $i\in\q$, is the truncated subsystem extracted from the balanced realization using the generalized structured Gramians. Then,
$$
\|\mathcal{F}(N,\hat{G}(s))-\mathcal{F}(N,G(s))\|_\infty\leq 2\sum_{i=1}^q \sum_{k=r_i+1}^{n_i}\sigma_{i,k}.
$$
\end{theorem}

This bound provides us with a simple procedure for finding an appropriate order for the reduced subsystem. First, we rank the Hankel singular values of each subsystem $\{\sigma_{i,j}\}_{j=1}^{n_i}$ in a descending order. Then, we may select the order of the reduced system by selecting the index after which there is a significant drop in the value of the singular values.

\subsection{Singular Perturbation}
Using the singular perturbation, the reduced subsystem $i\in\q$ is given by
$$
\hat{G}_i(s)=\hat{C}_i(sI-\hat{A}_i)^{-1}\hat{B}_i+\hat{D}_i,
$$
where 
\begin{align*}
\hat{A}_i&=\bar{A}_i^{11}-\bar{A}_i^{12}(\bar{A}_i^{22})^{-1}\bar{A}_i^{21},&
\hat{B}_i&=\bar{B}_i^{1}-\bar{A}_i^{12}(\bar{A}_i^{22})^{-1}\bar{B}_i^{2},\\
\hat{C}_i&=\bar{C}_i^{1}-\bar{C}_i^{2}(\bar{A}_i^{22})^{-1}\bar{A}_i^{21},&
\hat{D}_i&=\bar{D}_i-\bar{C}_i^{2}(\bar{A}_i^{22})^{-1}\bar{B}_i^{2}.
\end{align*}
A nice property of the singular perturbation is that the reduced system and the original system behave closely for low frequencies, that is, $\mathcal{F}(N,\hat{G}(0))=\mathcal{F}(N,G(0))$~\cite{green1995linear}.

Note that since, in an interconnected system, the other subsystems act as a low-pass filter for any given subsystem if the subsystems are all strictly proper (i.e., they have no direct term), matching the behavior of a subsystem at low frequencies might result in a better closed-loop performance (as the high frequencies are filtered out anyhow and, hence, match each other perfectly). Therefore, using singular perturbation in structured model reduction is, heuristically, better justified. 

\begin{figure}
\centering
\begin{tikzpicture}[scale=1,>=stealth']
\node (N) at (+0.0,+0.0) [draw,rectangle,fill=white,minimum width=2.0cm,minimum height=1.0cm] {$N$   };
\node (G) at (+0.0,-1.5) [draw,rectangle,fill=white,minimum width=2.0cm,minimum height=1.0cm] {$G(s)$};
\draw[->,line width=1pt] (+2.0,-0.2) to node[right]{} (+2.0,-1.5)
						 (+2.0,-1.5) to node[right]{} (+1.0,-1.5)
						 (+2.0,-0.2) to node[right]{} (+1.0,-0.2);
\draw[->,line width=1pt] (-1.0,-0.2) to node[right]{} (-2.0,-0.2)
						 (-2.0,-0.2) to node[right]{} (-2.0,-1.5)
						 (-2.0,-1.5) to node[right]{} (-1.0,-1.5);
\node (1) at (+2.5,+0.6) {$w(t)$};
\node (1) at (-2.5,+0.6) {$z(t)$};
\node (1) at (-2.5,-0.8) {$u(t)$};
\node (1) at (+2.5,-0.8) {$y(t)$};
\node (M) at (+0.0,-3.0) [draw,rectangle,fill=white,minimum width=2.0cm,minimum height=1.0cm] {$N$   };
\node (H) at (+0.0,-4.5) [draw,rectangle,fill=white,minimum width=2.0cm,minimum height=1.0cm] {$\frac{1}{s}I_r$};
\draw[->,line width=1pt] (+2.2,-3.2) to node[right]{} (+2.2,-6.7)
						 (+2.2,-6.7) to node[right]{} (+1.0,-6.7)
						 (+2.2,-3.2) to node[right]{} (+1.0,-3.2);
\draw[->,line width=1pt] (+2.0,-4.5) to node[right]{} (+2.0,-6.3)
						 (+2.0,-6.3) to node[right]{} (+1.0,-6.3)
						 (+2.0,-4.5) to node[right]{} (+1.0,-4.5);
\draw[->,line width=1pt] (-1.0,-3.2) to node[right]{} (-2.2,-3.2)
						 (-2.2,-3.2) to node[right]{} (-2.2,-6.7)
						 (-2.2,-6.7) to node[right]{} (-1.0,-6.7);
\draw[->,line width=1pt] (-1.0,-4.5) to node[right]{} (-2.0,-4.5)
						 (-2.0,-4.5) to node[right]{} (-2.0,-6.3)
						 (-2.0,-6.3) to node[right]{} (-1.0,-6.3);
\node (K) at (+0.0,-6.5) [draw,rectangle,fill=white,minimum width=2.0cm,minimum height=1.0cm] {$\Phi=\matrix{cc}{\hat{A}&\hat{B}\\ \hat{C}&\hat{D}}$};
\node (1) at (+2.5,-2.5) {$w(t)$};
\node (1) at (-1.8,-2.5) {$z'(t)$};
\node (1) at (-2.6,-2.5) {$+$};
\node (1) at (-3.4,-2.5) {$-$};
\node (S1) at (-3.0,-2.8) [draw,circle,fill=white,minimum width=0.1cm] {};
\draw[->,line width=1pt] (+3.0,+0.2) to node[right]{} (+3.0,-2.8)
						 (+4.0,-2.8) to node[right]{} (+1.0,-2.8);
\draw[->,line width=1pt] (-1.0,-2.8) to node[right]{} (S1);
\draw[->,line width=1pt] (S1) to node[right]{} (-4.5,-2.8);
\draw[->,line width=1pt] (+3.0,+0.2) to node[right]{} (+1.0,+0.2);
\draw[->,line width=1pt] (-1.0,+0.2) to node[right]{} (-2.0,+0.2)
						 (-1.0,+0.2) to node[right]{} (-3.0,+0.2)
						 (-3.0,+0.2) to node[right]{} (S1);
\draw[-,dashed,line width=1pt] (-3.7,+1.0) to node[right]{} (+3.4,+1.0)
							   (+3.4,+1.0) to node[right]{} (+3.4,-5.5)
							   (+3.4,-5.5) to node[right]{} (-3.7,-5.5)
							   (-3.7,-5.5) to node[right]{} (-3.7,+1.0);
\node (1) at (+0.0,+1.3) {$P(s)$};
\end{tikzpicture}
\caption{\label{fig:interconnection_Hinfityoptim} The realization of the transfer function $\mathcal{F}(N,G(s))-\mathcal{F}(N,\hat{G}(s))$ as a feedback interconnection with the unknowns in the feedback gain.\vspace{-.1in} }
\end{figure}
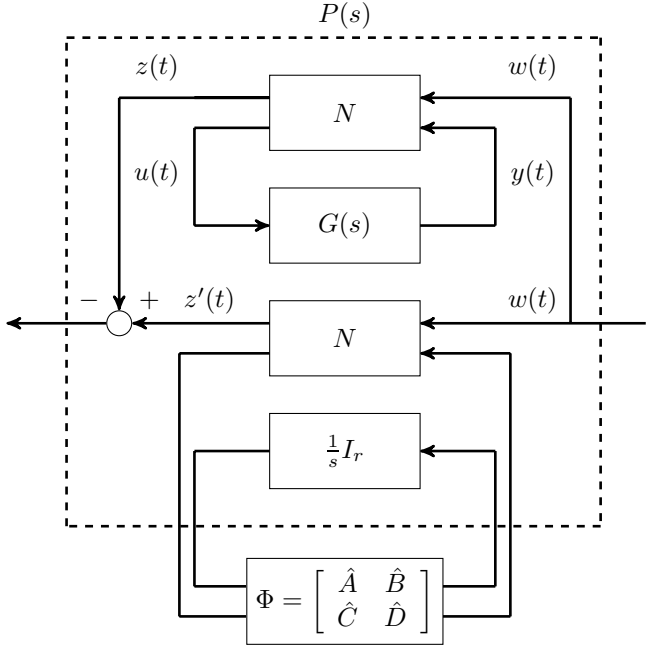

\section{Structured Model Reduction Using Subgradient Optimization} \label{sec:subgradient}
Model reduction using structured Gramians is a heuristic and, although very useful, it can fail occasionally. The standard balanced truncation algorithm (without the structure) has no $\mathcal{H}_\infty$ optimality property either. Moreover, model reduction methods using generalized structured Gramians can give solutions that are potentially far from the optimal solution. In addition, there is also no guarantee that, in general, the generalized structured Gramians even exist. Hence, in this section, we use subgradient optimization algorithm to improve the quality of the heuristic solutions. To do this, we use the methodology introduced in~\cite{1576856} to compute (sub)optimal $\mathcal{H}_\infty$ control laws using subgradient optimization algorithm. 

We can rewrite the error transfer function $\mathcal{F}(N,G(s))-\mathcal{F}(N,\hat{G}(s))$ as a feedback interconnection, where all the decision variables $(\hat{A},\hat{B},\hat{C},\hat{D})$, i.e., the model matrices of the reduced subsystems, are in the feedback gain. This is portrayed in Figure~\ref{fig:interconnection_Hinfityoptim}. Here, $r=r_1+\dots+r_q$ denotes the total order of the reduced subsystems. By definition (of the presented feedback interconnection), we have $\mathcal{F}(N,G(s))-\mathcal{F}(N,\hat{G}(s))=\mathcal{F}(P(s),\Phi)$. 
Now, following the results of~\cite{1576856}, we can easily construct the subdifferentials $\partial_{\Phi}\|\mathcal{F}(P(s),\Phi)\|_\infty$. To keep the matrices $\hat{A}$, $\hat{B}$, $\hat{C}$, and $\hat{D}$ block-diagonal (so as to preserve the subsystems and their interconnection structure), we should use the projected subgradients 
\begin{align*}
\Delta_i=\left\{\Xi\circ \matrix{cc}{\Psi_A & \Psi_B \\ \Psi_C & \Psi_D}\bigg|\, \Xi\in\partial_{\Phi}\|\mathcal{F}(P(s),\Phi)\|_\infty \right\},
\end{align*}
where 
\begin{align*}
\Psi_A&=\diag(\mathbf{1}_{r_1\times r_1},\dots,\mathbf{1}_{r_q\times r_q}),\\
\Psi_B&=\diag(\mathbf{1}_{r_1\times m_1},\dots,\mathbf{1}_{r_q\times m_q}),\\
\Psi_C&=\diag(\mathbf{1}_{p_1\times r_1},\dots,\mathbf{1}_{p_q\times r_q}),\\
\Psi_D&=\diag(\mathbf{1}_{p_1\times m_1},\dots,\mathbf{1}_{p_q\times m_q}).
\end{align*}
Here, $X\circ Y$ denotes the Hadamard product, also known as the element-wise product, of matrices $X$ and $Y$ with appropriate dimensions. Now, we can propose a numerical algorithm to construct a locally optimal reduced system by moving in the opposite direction of this projected subgradients. Such a numerical algorithm is discussed in length in Section~VI.F in~\cite{1576856}.

%the entrywise produce of the subgradient matrix $\partial_{\Phi}\|\mathcal{F}(P(s),\Phi)\|_\infty$ and the matrices $\Psi_A$, $\Psi_B$, $\Psi_C$, and $\Psi_D$ in Line~4 of Algorithm~\ref{alg:1}.

%\begin{algorithm}
%\caption{\label{alg:1} Model reduction algorithm using the subgradient optimization method.}
%\begin{algorithmic}[1]
%%\REQUIRE $\epsilon\in\mathbb{R}_{>0}$, $\{\mu_i\}_{i\in\mathbb{N}_0}$;
%%\ENSURE The reduced subsystems;
%\STATE Initialize $(\hat{A}^0,\hat{B}^0,\hat{C}^0,\hat{D}^0)$;
%\FOR{$i=1,2,\dots$}
%\STATE{ Calculate 
%\begin{align*}
%\hspace{-.1in}\Delta_i=\partial_{\Phi}\|\mathcal{F}(P(s),\Phi)\|_\infty\big|_{\Phi=\scriptsize\matrix{cc}{\hat{A}^{i}&\hat{B}^{i}\\ \hat{C}^{i}&\hat{D}^{i}}}\circ \matrix{cc}{\Psi_A & \Psi_B \\ \Psi_C & \Psi_D}\hspace{-.05in},
%\end{align*}
%where 
%\begin{align*}
%\Psi_A&=\diag(\mathbf{1}_{r_1\times r_1},\dots,\mathbf{1}_{r_q\times r_q}),\\
%\Psi_B&=\diag(\mathbf{1}_{r_1\times m_1},\dots,\mathbf{1}_{r_q\times m_q}),\\
%\Psi_C&=\diag(\mathbf{1}_{p_1\times r_1},\dots,\mathbf{1}_{p_q\times r_q}),\\
%\Psi_D&=\diag(\mathbf{1}_{p_1\times m_1},\dots,\mathbf{1}_{p_q\times m_q});
%\end{align*}
%}
%\STATE{ Update
%\begin{align*}
%\matrix{cc}{\hat{A}^{i+1}&\hat{B}^{i+1}\\ \hat{C}^{i+1}&\hat{D}^{i+1}}
%=&\matrix{cc}{\hat{A}^{i}&\hat{B}^{i}\\ \hat{C}^{i}&\hat{D}^{i}}-\mu_i \Delta_i;
%\end{align*}
%}
%\ENDFOR
%\end{algorithmic}
%\end{algorithm}
%
%\begin{theorem} 
%
%\end{theorem}
%
%\begin{proof} 
%
%\end{proof}

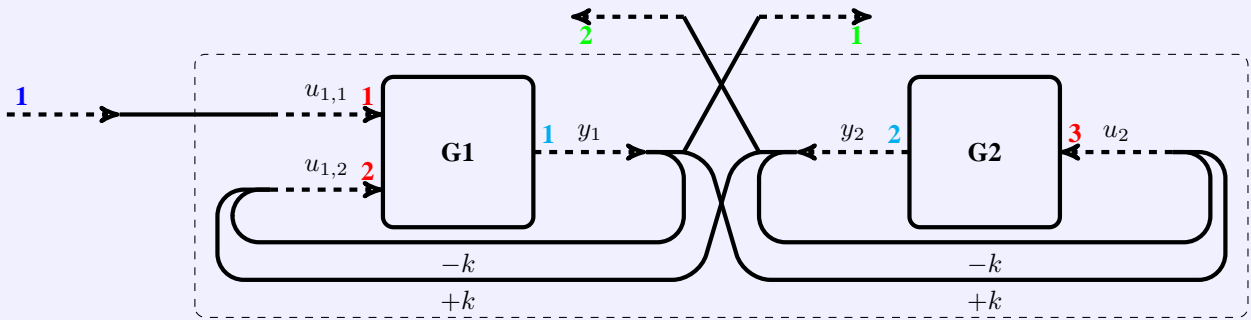
\begin{figure*}[t]
\centering
\begin{tikzpicture}[>=stealth']
\node (N) at (+0.0,+0.0) [draw,rectangle,fill=blue!6,rounded corners=10pt] {
\begin{minipage}{0.99\textwidth}
{\large \textbf{Structured ModeL reductIon (SiMpLIfy) Toolbox}} \\[0.4em]
SiMpLIfy is a MATLAB toolbox for structured model reduction of interconnected systems. The user guide and the m-files for this toolbox can be found in {\small\url{http://simplifytoolbox.tumblr.com/}}. There are also several demos, including the numerical example in this paper, attached to the toolbox files. In what follows, we briefly introduce the syntax for using this toolbox in the context of the mass--spring mechanical example utilized in this paper.

The first step is to construct an interconnected system as an instance of a `SystemNetwork' class defined in the toolbox as
\begin{verbatim}
>> iedges=[1 2 −k; 2 2 k; 1 3 k; 2 3 −k];
>> einedges=[1 1];
>> eoutedges=[1 1; 2 2];
>> eedges=[ ];
>> systemNetwork=SystemNetwork(iedges,einedges,eoutedges,eedges,G1,G2)
\end{verbatim}
In this example, `iedges', `einedges', `eoutedges', and `eedges' list, repsectively, the edges between interal  outputs and internal inputs, the edges between external inputs and internal inputs, the edges between internal ouputs and external outputs, and the edges between external inputs and external outputs. See Figure~\ref{fig:matlab:networkStructure} for a visualization of these edges. Moreover, `G1' and `G2' are the transfer functions of the first and the second subsystems. Now, we may use the following command to visualize the Hankel singular values and to compare them with regular Hankel singular values
\begin{verbatim}
>> compareHankels(systemNetwork)
\end{verbatim}
To reduce the order of the subsystems, we can use balanced reduction with structured Gramians as
\begin{verbatim}
>> red=balancedNetworkReduction(systemNetwork,[6 3])
\end{verbatim}
We may also use the following command to extract the reduced subsystems using singular perturbation
\begin{verbatim}
red = balancedNetworkReduction(systemNetwork,[6 3],`ReductionMethod',`perturbation')
\end{verbatim}
This command construct subsystems of orders 6 and 3, repsectively. After the reduction, the subsystems can be extracted using the commands
\begin{verbatim}
>> red.extractSubsystem(1)
>> red.extractSubsystem(2)
\end{verbatim}
Finally, we can improve the quality of the reduced models by using the subgradient optimization algorithm as
\begin{verbatim}
>> optred=improveNetworkReduction(systemNetwork,red)
\end{verbatim}
This command supplies the preliminary reduced model `red', extracted from the balanced truncation, as an intial point to the numerical algorithm.  We invite the interested readers to check the user guide developed for the toolbox to learn about all the other implemented algorithms. We have also included `testMechSystem.m' in the toolbox for creating this mechanical system and reducing it.
\vspace{.1in}
\begin{center}
\begin{tikzpicture}
\draw [black,dashed,rounded corners] (-2.5,2.8) rectangle (11.5,6.3);
\draw [black,ultra thick,rounded corners] (0,4) rectangle (2,6);
\draw [black,ultra thick,rounded corners] (7,4) rectangle (9,6);
\node at (1,5) {\bf G1};
\node at (8,5) {\bf G2};
\node [blue ,above] at (-4.8,+5.5) {\bf 1};
\node [red  ,above] at (-0.2,+5.5) {\bf 1};
\node [red  ,above] at (-0.2,+4.5) {\bf 2};
\node [red  ,above] at (+9.2,+5.0) {\bf 3};
\node [green,below] at (+2.7,+6.8) {\bf 2};
\node [green,below] at (+6.3,+6.8) {\bf 1};
\node [cyan ,above] at (+6.8,+5.0) {\bf 2};
\node [cyan ,above] at (+2.2,+5.0) {\bf 1};
\draw [->,dashed,ultra thick] (-5.0 ,+5.5) -- (-3.5,+5.5);
\draw [->,dashed,ultra thick] (+5.0 ,+6.8) -- (+6.5,+6.8);
\draw [->,dashed,ultra thick] (+4.0 ,+6.8) -- (+2.5,+6.8);
\draw [->,dashed,ultra thick] (-1.5 ,+5.5) -- node[above] {$u_{1,1}$} (+0.0,+5.5);
\draw [->,dashed,ultra thick] (-1.5 ,+4.5) -- node[above] {$u_{1,2}$} (+0.0,+4.5);
\draw [->,dashed,ultra thick] (+10.5,+5.0) -- node[above] {$u_{2}$  } (+9.0,+5.0);
\draw [->,dashed,ultra thick] (+7.0 ,+5.0) -- node[above] {$y_{2}$  } (+5.5,+5.0);
\draw [->,dashed,ultra thick] (+2.0 ,+5.0) -- node[above] {$y_{1}$  } (+3.5,+5.0);
\draw [- ,ultra thick] (-3.5,+5.5) -- (-1.5,+5.5);
\draw [- ,ultra thick] (+3.5,+5.0) -- (+4.0,+5.0) -- (+4.0,+3.8) -- node[below] {$-k$} (-2.0,+3.8) -- (-2.0,+4.5) -- (-1.5,+4.5);
\draw [- ,ultra thick] (+5.5,+5.0) -- (+5.0,+5.0) -- (+5.0,+3.8) -- node[below] {$-k$} (11.0,+3.8) -- (11.0,+5.0) -- (10.5,+5.0);
\draw [- ,ultra thick] (+3.5,+5.0) -- (+4.3,+5.0) -- (+4.8,+3.3) -- node[below] {$+k$} (11.2,+3.3) -- (11.2,+5.0) -- (10.5,+5.0);
\draw [- ,ultra thick] (+5.5,+5.0) -- (+4.7,+5.0) -- (+4.2,+3.3) -- node[below] {$+k$} (-2.2,+3.3) -- (-2.2,+4.5) -- (-1.5,+4.5);
\draw [- ,ultra thick] (+4.0,+6.8) -- (+5.0,+5.0);
\draw [- ,ultra thick] (+5.0,+6.8) -- (+4.0,+5.0);
\end{tikzpicture}
\caption{The network structure for the mechanical system. The nodes \textbf{G1} and \textbf{G2} denote the subsystems. The internal inputs (the inputs of each subsystem) are numbered sequentially in {\color{red}red} color. The internal outputs (the outputs of each subsystem) are numbered sequentially in {\color{cyan}cyan} color. The external inputs (the inputs of the overall interconnected system) are numbered sequentially in {\color{blue}blue} color. Finally, the external outputs (the outputs of the overall interconnected system) are numbered sequentially in {\color{green}green} color. The solid curves portray the edges between internal inputs, internal outputs, external inputs, and external outputs. The weights on the edges are displayed only if they are not equal to identity. Notice that if more than one edge is going to an input, implicitly, we mean that these edges are summed together before being fed to that input.  \label{fig:matlab:networkStructure}}
\end{center}
\end{minipage}
};
\end{tikzpicture}
\vspace{-.3in}
\end{figure*}

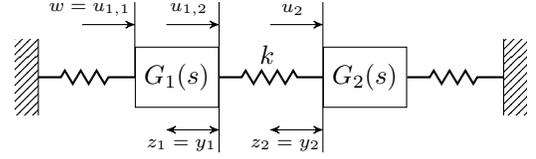
\begin{figure}[t!]
\centering
\begin{tikzpicture}[>=stealth']
\tikzstyle{spring}=[thick,decorate,decoration={zigzag,pre length=0.3cm,post length=0.3cm,segment length=6}]
\tikzstyle{ground}=[fill,pattern=north east lines,draw=none,minimum width=0.75cm,minimum height=0.3cm]
\node (wall1) [ground, rotate=-90, minimum width=1cm,yshift=-2cm] {};
\draw (wall1.north east) -- (wall1.north west);
\node (M1) [rectangle,draw,minimum width=1cm, minimum height=0.8cm] {$G_1(s)$};
\draw [spring] (wall1) -- (M1);
\node (M2) [rectangle,draw,minimum width=1cm, minimum height=0.8cm] at (2.5,0) {$G_2(s)$};
\draw [spring] (M1) -- (M2);
\node (wall2) [ground, rotate=-90, minimum width=1cm,yshift=4.5cm] {};
\draw (wall2.south east) -- (wall2.south west);
\draw [spring] (M2) -- (wall2);
\node [] at (1.2,.3) {$k$};
\draw (-0.56,0) -- (-0.56,1);
\draw (+0.56,0) -- (+0.56,1);
\draw (+0.56,0) -- (+0.56,-1);
\draw (+1.94,0) -- (+1.94,1);
\draw (+1.94,0) -- (+1.94,-1);
\draw [->] (-1.26,0.7) -- (-0.56,0.7);
\node [] at (-1.16,0.9) {\scriptsize $w=u_{1,1}$};
\draw [->] (-0.14,0.7) -- (+0.56,0.7);
\node [] at (0.16,0.9) {\scriptsize $u_{1,2}$};
\draw [<->] (-0.14,-0.7) -- (+0.56,-0.7);
\node [] at (0.06,-0.9) {\scriptsize $z_1=y_1$};
\draw [->] (1.24,0.7) -- (+1.94,0.7);
\node [] at (1.54,0.9) {\scriptsize $u_{2}$};
\draw [<->] (1.24,-0.7) -- (+1.94,-0.7);
\node [] at (1.44,-0.9) {\scriptsize $z_2=y_2$};
\end{tikzpicture}
\caption{ A mass--spring mechanical system. \label{fig:mechexample}}
\end{figure}

\section{Numerical Example} \label{sec:numerical}
\begin{figure}[t!]
\centering
\includegraphics[width=0.95\columnwidth]{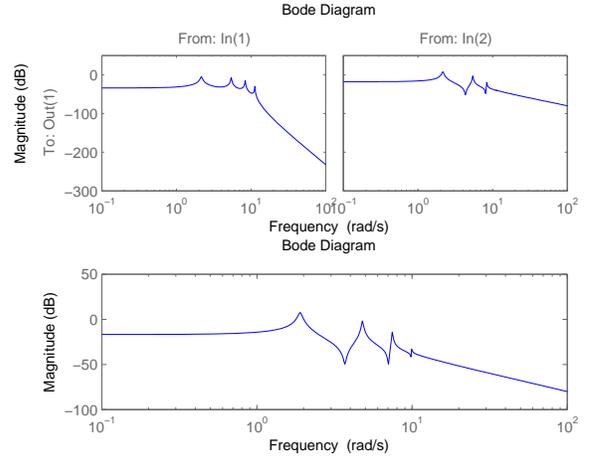}
\caption{ Bode magnitude plots of individual subsystems $G_1(s)$ (top) and $G_2(s)$ (bottom).\label{fig:bode}
\vspace{-.2in}}
\end{figure}

Let two elastic masses be interconnected by a linear spring with constant $k=10$ as in Figure~\ref{fig:mechexample}. The first mass model, $G_1(s)$, takes the forces $u_{1,1}(t)$ and $u_{1,2}(t)$ as inputs and provides its positions as an output $y_1(t)$. Similarly, the second mass model, $G_2(s)$, takes the force $u_{2}(t)$ as an input and provides the its positions as an output $y_2(t)$.  These models are of high degree as they are discretized partial differential equations of elastic bodies, i.e., $\deg(G_1(s))=8$ and $\deg(G_2(s))=10$. The Bode magnitude plots of these systems\footnote{The m-file for constructing these models is included in the toolbox as a demo.} are shown in Figure~\ref{fig:bode}.

The models are interconnected with a linear spring, and the mapping of interest (the transfer function from the external input to the external output) could be how the force $w(t)=u_{1,1}(t)$ maps to the two positions $z(t)=\matrix{cc}{y_1 (t) & y_2 (t)}^\top$. The network is, hence, modeled as
\begin{align*}
N
&=\matrix{c:c}{D_E & D_F \\ \hdashline D_H & D_K}=\matrix{c:cc}{0 & 1 & 0 \\ 0 & 0 & 1 \\  \hdashline 1 & 0 & 0 \\ 0 & -k & k \\ 0 & k & -k}.
\end{align*}

\begin{figure}[t!]
\centering
\includegraphics[width=0.85\columnwidth]{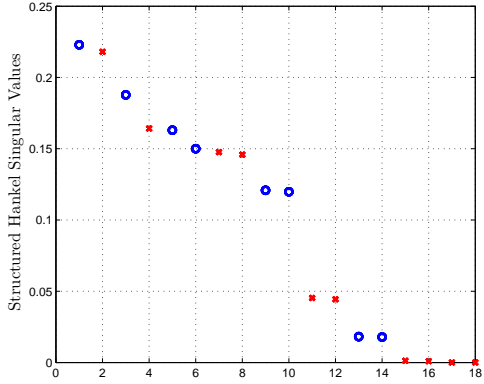}
\caption{The structured Hankel singular values for the 18th-order interconnected model $\mathcal{F}(N,G(s))$. 
The singular values corresponding to $G_1(s)$ and $G_2(s)$ are, respectively, marked with `${\color{blue}\circ}$' and `${\color{red}\times}$'. \label{fig:structuredHSV} \vspace{-.1in}}
\end{figure}

\begin{table}
\caption{\label{table:1} The reduction error for the balanced truncation with the structured Gramians.}
\scriptsize
\centering
\begin{tabular}{|c||c|c|c|c|c|}
\hline 
$r_2/r_1$ & 8 & 6 & 4 & 2 \\ \hline \hline 
10 		  & $0.0$ & $3.58\times 10^{-2}$ & $2.40\times 10^{-1}$ & $\infty$ \\ \hline 
8  		  & $1.23\times 10^{-5}$ & $3.58\times 10^{-2}$ & $2.40\times 10^{-1}$ & $\infty$ \\ \hline 
6  		  & $7.62\times 10^{-3}$ & $3.58\times 10^{-2}$ & $2.40\times 10^{-1}$ & $\infty$ \\ \hline 
4  		  & $2.08\times 10^{-1}$ & $2.07\times 10^{-1}$ & $3.28\times 10^{-1}$ & $\infty$ \\ \hline 
2  		  & $3.98\times 10^{-1}$ & $3.98\times 10^{-1}$ & $4.20\times 10^{-1}$ & $\infty$ \\ \hline 
\end{tabular}
\caption{\label{table:2} The reduction error for the model reduction method using the subgradient optimization algorithm.}
\begin{tabular}{|c||c|c|c|c|c|}
\hline 
$r_2/r_1$ & 8 & 6 & 4 & 2 \\ \hline \hline 
10 		  & $0.0$ & $2.17\times 10^{-3}$ & $6.33\times 10^{-3}$ & $1.23\times 10^{-1}$ \\ \hline 
8  		  & $6.33\times 10^{-7}$ & $4.57\times 10^{-3}$ & $3.10\times 10^{-2}$ & $1.21\times 10^{-1}$ \\ \hline 
6  		  & $9.99\times 10^{-4}$ & $3.08\times 10^{-2}$ & $5.66\times 10^{-2}$ & $1.37\times 10^{-1}$ \\ \hline 
4  		  & $3.42\times 10^{-2}$ & $5.68\times 10^{-2}$ & $1.40\times 10^{-1}$ & $1.79\times 10^{-1}$ \\ \hline 
2  		  & $1.36\times 10^{-1}$ & $1.38\times 10^{-1}$ & $1.80\times 10^{-1}$ & $1.86\times 10^{-1}$ \\ \hline 
\end{tabular}
\vspace{-.15in}
\end{table}

\begin{figure}[t!]
\centering
\includegraphics[width=0.85\columnwidth]{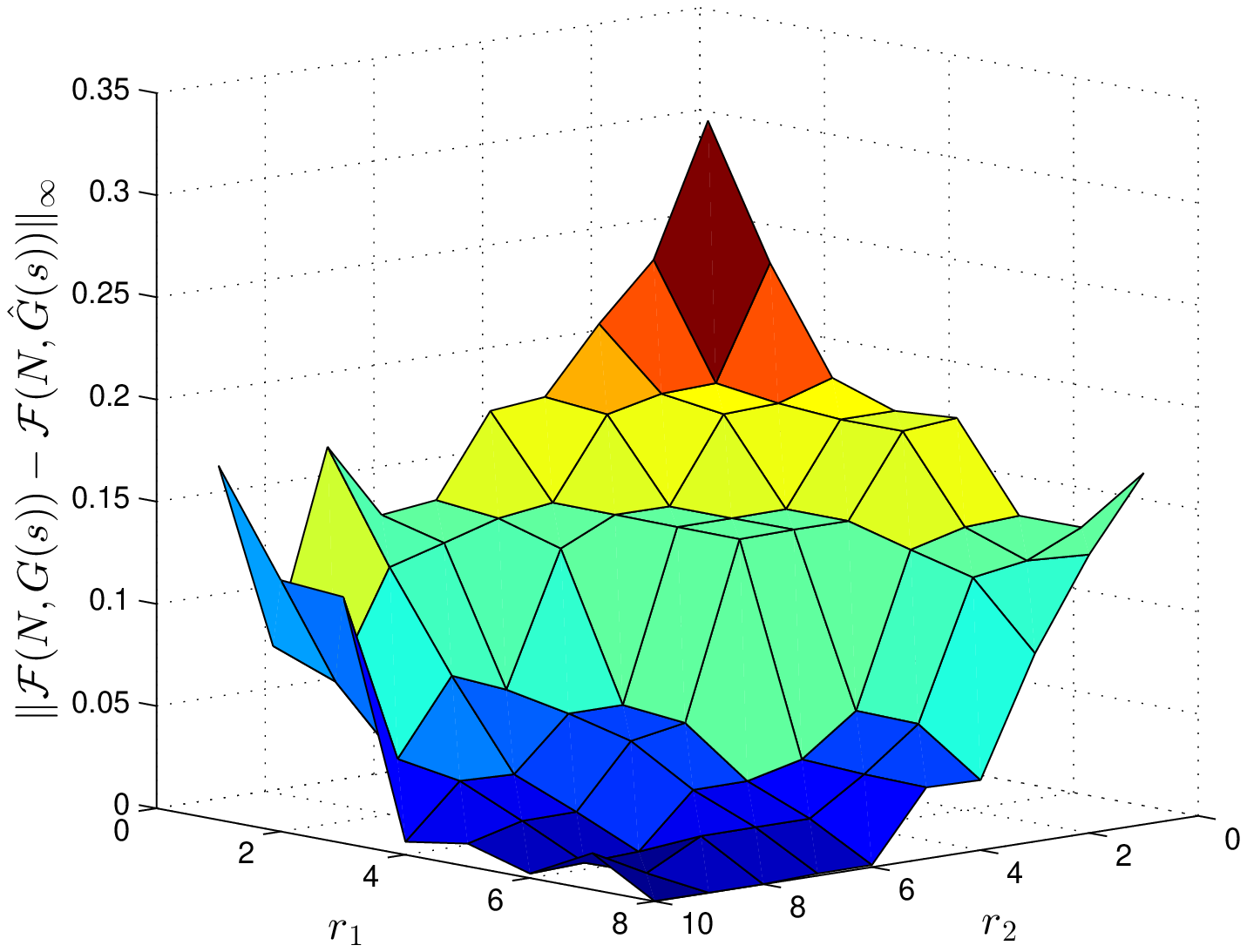}
\caption{ Reduction error $\|\mathcal{F}(N,G(s))-\mathcal{F}(N,\hat{G}(s))\|_\infty$ for the model reduction method using the subgradient optimization algorithm as function of the reduced systems' orders.  \label{fig:error_subgradient}}
\vspace{-.2in}
\end{figure}

Figure~\ref{fig:structuredHSV} illustrates the structured Hankel singular values extracted from the structured Gramians. The singular values corresponding to the first subsystem, $\{\sigma_{1,k}\}_{k=1}^8$, are marked with `${\color{blue}\circ}$' and the singular values corresponding to the second subsystem, $\{\sigma_{2,k}\}_{k=1}^{10}$, are marked with `${\color{red}\times}$'. For instance, these structured singular values show that the order of the second subsystem can be reduced to 6 without sacrificing the performance (i.e., preserving a similar input-output behavior).

Table~\ref{table:1} shows the model reduction error $\|\mathcal{F}(N,G(s))-\mathcal{F}(N,\hat{G}(s))\|_\infty$ for various selections of $r_1$ and $r_2$, i.e., the order of the reduced subsystems, when using the balanced truncation with structured Gramians. We can evidently see that reducing the order of the second subsystem to 6 does not introduce much error, which certifies our intuition from the structured Hankel singular values. Because the structured Gramians give heuristic reduction methods, we cannot expect a stable interconnected system with the reduced subsystems. This is evident from the last column of Table~\ref{table:1}.

To improve the quality of the reduced models, we can use the model reduction methods using subgradient optimization algorithm. We initialize this numerical algorithms with the reduced model from the balanced truncation with structured Gramians. Table~\ref{table:2} shows the reduction error $\|\mathcal{F}(N,G(s))-\mathcal{F}(N,\hat{G}(s))\|_\infty$ for various selections of $r_1$ and $r_2$ in this case. We can easily see that the results of this algorithm are much better than the ones extracted from the balanced truncation with structured Gramians. The reduction is error is also portrayed in Figure~\ref{fig:error_subgradient} for various orders. If we were to recover the global optimum (with the optimization algorithm), the reduction error would have been a decreasing function of the order. However, although the error is mostly decreasing with increasing the orders, this not true for all cases, which is because the proposed algorithm at best recovers a local optimum.

\section{Conclusions} \label{sec:conclusions}
We presented a toolbox for structured model reduction in MATLAB. It contains model reduction algorithms based on balanced truncation and singular perturbation. To construct the balanced realization of the subsystems, we use structured Gramians and generalized structured Gramians. The latter resulted in bounds on the reduction error. We also proposed a model reduction algorithm using a subgradient optimization algorithm. The algorithms were compared on a structured mechanical system. Future research can focus on extending the modules in the toolbox to admit uncertain models and/or
parameter-dependent ones.

\bibliographystyle{ieeetr}
\bibliography{ref}
\end{document}